\documentclass[12pt]{article}
\usepackage{graphicx}
\usepackage{amssymb}
\usepackage{epstopdf}
\title{Corrig\' e}
\author{}
   
\usepackage{graphicx}
\usepackage{amssymb}
\usepackage{epstopdf}
\title {A family of two dimensionally determined ergodic processes }
\date{}
\author {Doureid Hamdan} 
\title
  {Ergodic quasi-invariant measures on  topologically mixing subshifts     are isomorphic to  Bernoulli shifts}
\begin{document}
\maketitle

{\abstract{\textwidth=4,5 in} We prove that a shift ergodic measure on a topologically mixing sub-shift  is isomorphic to a Bernoulli shift whenever it is quasi invariant under permutations of finite number of coordinates. We prove also that Gibbs measures on topologically mixing subshift of finite type are quasi invariant.}


\section{Introduction, Notation}
Usually a sequence    $\hskip 0,1 cm (X_n)_{n\in \mathbb Z}\hskip 0,1 cm $ of random variables  is said exchangeable if   the law $P_\sigma$ of the process $(X_{\sigma(n)})_{n\in \mathbb Z}$ is equal  the law $P$ of the process $(X_n)_{n\in \mathbb Z}$, for every permutation $\sigma$  of the set $\mathbb Z$ of integers, which leaves fixed all but a finite number of integers. Naturally, $\hskip 0,1 cm (X_n)_{n\in \mathbb Z}\hskip 0,1 cm $ will then be said quasi exhangeable if $P_\sigma$ is equivalent to $P$, for any such $\sigma$.\\ 
In $[6]$ Hewitt and Savage enlarge the category of the  state space and obtain a generalization of De Finetti's Theorem which says that an an exchangeable sequence of random variables is a mixture of  i.i.d. sequences. We note that exchangeability implies stationarity. In [5] it is proved that, if the dynamical system generated by a stationary  quasi exchangeable process is ergodic, then it is isomorphic to a Bernoulli process, and in the particular case where the family of all Radon-Nikodym derivatives is uniformly integrable, then the process is exchangeabe. Also an  application of one result from $[5]$ yields  that the translation invariant  determinantal processes considered in [2] are isomorphic to Bernoulli processes.\\  In this paper we shall   consider only finite state space stationary processes
$(X_n)_{n\in \mathbb Z}$ whose underlying  dynamical system is a topologically mixing subshift.   We prove that  they are isomorphic to Bernoulli shifts, whenever  they are quasi exchangeable and ergodic, a result which generalizes Theorem 1 of $[5]$.\\
First we establish the terminology and   notation which will be needed later. \\We suppose that,  for every $n\in \mathbb Z$,  $X_n$ is the $n^{\textit{th}}$ coordinate on the product space $\Omega:=K^\mathbb Z$, where $K$ is the finite state space, and   also that the law $P$ is concentrated on a subshift in $\Omega$. \\ 
 Then  
\begin{eqnarray*}
{\cal F}:= \textit{the sigma algebra of all subsets of}\hskip 0,2 cm K,\\
\Omega:=K^\mathbb Z,\hskip 0,1 cm \textit{ is endowed with the product sigma algebra} \hskip 0,3 cm {\cal E}:={\cal F}^{\otimes \mathbb Z},\hskip 0,2 cm  \textit{  and}\\\hskip 0,2 cm       X_n(\omega)=\omega_n,\hskip 0,2 cm  \textit{ for all integer }\hskip 0,05 cm n\in \mathbb Z\hskip 0,3 cm \textit{and all}\hskip 0,4 cm \omega\in \Omega. \end{eqnarray*}
    The shift transformation is denoted  $\hskip 0,1 cm S:\hskip 0,1 cm (S\omega)_n=\omega_{n+1},\hskip 0,1 cm n\in  \mathbb Z,\hskip 0,1 cm \omega \in \Omega. $ \\
Let   $G$ be the group of all permutations of $\hskip 0,1 cm \mathbb Z\hskip 0,1 cm $  and
 $\hskip 0,1 cm H\subset G,$ be  the subgroup of all permutations  with finite support:  \begin{eqnarray*}
\sigma\in H\iff \sigma\in G, \hskip 0,2 cm \textit{and} \hskip 0,3 cm  \exists N, \hskip 0,2 cm  \sigma(n)=n, \forall n,  \mid n\mid \ge N.\end{eqnarray*}
For any $\tau\in G, \hskip 0,2 cm$ let  the transformation $T_\tau:\Omega\rightarrow \Omega$, be defined for all  $ \hskip 0,2 cm \omega\in \Omega, \hskip 0,2 cm$  by 
 \begin{eqnarray}(T_\tau(\omega))_n=\omega_{\tau(n)}, \hskip 0,2 cm \forall  n\in \mathbb Z.\end{eqnarray} 
 Then $T_\tau$ is ${\cal E}$-measurable, and, when $K$ is a topological space,  $T_\tau$ is continuous for the product topology on  $\Omega$.  Also     for all $\sigma$ and $\tau$ in $H$,   \begin{eqnarray}T_{\tau\circ \sigma}=T_\sigma\circ T_\tau,\end{eqnarray} 
 from which follows that \begin{eqnarray}T_\sigma^{-1}=T_{\sigma^{-1}}.\end{eqnarray}
   Recall that then the quasi exchangeability of the sequence $(X_n)_{n\in \mathbb Z}$ is equivalent to the following
\\
{\bf{Definition 1}}\\
\textit{We say that a sequence  $(X_n)_{n\in \mathbb Z}$ of random variables, with law $\mu $,  is quasi exchangeable if  $\mu\circ T_\sigma^{-1}$ is equivalent to $\mu$, for all permutation $\sigma\in H$.}\\

In this case we denote the Radon-Nikodym derivative of $ \hskip 0,1 cm \mu\circ T_\sigma^{-1} \hskip 0,1 cm$ with respect to $ \hskip 0,2 cm \mu, \hskip 0,1 cm$ by  $ \hskip 0,1 cm \phi_\sigma$
\begin{eqnarray}
\phi_\sigma:={d\mu\circ T_\sigma^{-1}\over d\mu}.
\end{eqnarray}

 Before setting the precise statement, we establish some more notation. \\
  We endow  $ \Omega$  with the product topology of the discrete topologies on $K$.  Let $\Sigma$ be a closed shift-invariant subset of $\Omega$. We assume that $\forall j\in K$, there exists $x\in \Sigma$ such that $x_0=j$. \\
 We assume also that the system $(\Sigma,S)$ is topologically mixing, which means that  \textit{  for non empty open sets $\hskip 0,1 cm U,V\hskip 0,1 cm $ of $\Sigma$, there is an $\hskip 0,1 cm N\hskip 0,1 cm $ such that $\hskip 0,1 cm U\cap S^mV\ne\emptyset,\hskip 0,1 cm$ for all $\hskip 0,1 cm m\ge N$.}\\
 We shall also use the following notations.\\
If  $L$ and $s$ are integers with $L\ge 0$ and $s\ge 1$,  and   if $A_{-L},...,A_s$ are  measurable subsets  of $K$, we set 
\begin{eqnarray}
{\Pi}(A_{-L},...,A_0):=\{\omega\in \Sigma: \omega_{-L}\in A_{-L},...,\omega_0\in A_0\}\\
F(A_1,...,A_s):=\{\omega \in \Sigma:\omega_1\in A_1,...,\omega_s\in A_s\} ,
\end{eqnarray}
and for all $I\subset \mathbb Z$,
\begin{eqnarray*}
{\cal A}_I:= \textit{the smallest algebra containing the  sets} \hskip 0,2 cm \{\omega\in \Sigma:\omega_j\in A\}, j\in I, A\in {\cal F},\\
\textit{and}\hskip 2 cm 
{\cal B}_I:=\textit{the sigma algebra generated by }\hskip 0,2 cm  {\cal A}_I.
\end{eqnarray*}
In the following particular cases:
\begin{eqnarray*}
I=\{n\in \mathbb Z:n\le 0\}\hskip 0,4 cm  {\cal A}_I \hskip 0,2 cm  \textit{ is denoted}\hskip 0,2 cm  {\cal A}_{\le 0}\\
I=\{n\in \mathbb Z:n\ge  p\}\hskip 0,4 cm  {\cal A}_I \hskip 0,2 cm  \textit{ is denoted}\hskip 0,2 cm  {\cal A}_{\ge p}\\
I=\mathbb Z \hskip 0,4  cm  {\cal A}_I \hskip 0,2 cm  \textit{ is denoted}\hskip 0,2 cm  {\cal A}.
\end{eqnarray*}
The same notation will be used for ${\cal B}_I$, in particular  $\hskip 0,3 cm {\cal B}_{\mathbb Z}={\cal B}$. \\
Also, the smallest sigma-algebra which contains  two sigma-algebras   ${\cal F}_1$ and ${\cal F}_2$  is denoted by $\hskip 0,2 cm {\cal F}_1\vee {\cal F}_2$,  and the smallest one containing a family $\{{\cal F}_j:j\in J\}$, of sigma-algebras is denoted $\bigvee_{j\in J} {\cal F}_j$. \\
 For any  $\sigma\in H$, let us consider the subset  \begin{eqnarray}
\Gamma_\sigma:=\Sigma\cap T_\sigma^{-1}\Sigma.
\end{eqnarray}
  Then   $\Gamma_\sigma$ is a clopen set in $\Sigma$ and the following implication holds
   \begin{eqnarray*}
 \sigma \hskip 0,1 cm \textit{is an involution} \hskip 0,1 cm \Rightarrow  T_\sigma \Gamma_\sigma=\Gamma_\sigma=T_\sigma^{-1}\Gamma_\sigma.
\end{eqnarray*}

Let \begin{eqnarray}
H_+:=\{\sigma\in H:  \sigma \hskip 0,1 cm \textit{is an involution, } \hskip 0,1 cm   \Gamma_\sigma\ne \emptyset\},
\end{eqnarray}

 and,  for every $a\in K$, 
\begin{eqnarray}
\nonumber L(a)=\{j\in K:\{x\in \Sigma: x_0=j,x_1=a\}\ne \emptyset\}\hskip 1 cm \\\textit{and}\hskip 1 cm R(a)=\{j\in K:\{x\in \Sigma: x_0=a,x_1=j\}\ne \emptyset \}\}.
\end{eqnarray}
[ so that  $\hskip 0,1 cm \hskip 0,1 cm j\in L(a)\iff a\in R(j)$]. \\
Let \begin{eqnarray*}
{\cal Q}:=\{(a,b)\in K\times K:a\in L(b)\}
\end{eqnarray*}
 Let us also use the notation
\begin{eqnarray}
(x_j\in A_j, j\in J):=\{x\in \Sigma:x_j\in A_j, j\in J\}, \hskip 0,5 cm J\subset \mathbb Z,
\end{eqnarray}
where, for every $j$,  $A_j$ is a subset of $K,\hskip 0,1 cm$and in particular,  we set
\begin{eqnarray}
C_k(a,b,c,d)
:=(x_0=a,x_1=b,x_k=c,x_{k+1}=d)
\end{eqnarray}
\begin{eqnarray}
D_k(a,b,c,d):
=(x_0\in Lb,x_1\in Ra,x_k\in Ld,x_{k+1}\in Rc)
\end{eqnarray}
so that, with 
 \begin{eqnarray}
A=A(a,b):=(x_0=a,x_1=b),\hskip 1 cm \tilde{A}(a,b):=(x_0\in Lb, x_1\in Ra)\\
B=A(c,d):=(x_0=c,x_1=d),\hskip 1 cm \tilde{B}:=(x_0\in Ld,x_1\in Rc).
\end{eqnarray}
it is  clear that for $(a,b),(c,d)\in {\cal Q}$, we have the inclusions
\begin{eqnarray}
A(a,b)\subset \tilde{A}(a,b),
\end{eqnarray}
and also the following equalities 
 \begin{eqnarray}
C_k=C_k(a,b,c,d)=A\cap S^{-k}B, \hskip 1 cm 
D_k=D_k(a,b,c,d)=\tilde{A}\cap S^{-k}\tilde{B}.
\end{eqnarray}
Hence
\begin{eqnarray}
C_k(a,b,c,d)\subset D_k(a,b,c,d).
\end{eqnarray}
By topological mixing and because ${\cal Q}$ is finite, there exists $k_0$ such that 
 \begin{eqnarray}
C_k(a,b,c,d)\ne \emptyset,\hskip 0,2 cm \forall k\ge k_0,  \hskip 0,2 cm \forall (a,b),(c,d)\in {\cal Q}.
\end{eqnarray}
Hence, also, there exists $P_0(k)$ such that
 \begin{eqnarray}
C_k(a,b,c,d)\cap S^{-P}D_k(a,b,c,d)\ne \emptyset,\hskip 0,2 cm \forall k\ge k_0,  \hskip 0,2 cm \forall (a,b),(c,d)\in {\cal Q}, \forall P\ge P_0(k).
\end{eqnarray}
 
For all natural numbers  $ k$ and $ P,\hskip 0,2 cm $ such that $\hskip 0,2 cm  1\le k<P,\hskip 0,2 cm $ let us consider the permutation ( involution) $\hskip 0,2 cm  \sigma:=\sigma_{P,k}\in H\hskip 0,2 cm  $   which translates the "interval" $\hskip 0,1 cm I:=\mathbb N\cap [1,...,k]\hskip 0,1 cm $  by $\hskip 0,1 cm P$,  translates the $^"$interval$^"$ $\hskip 0,1 cm P+I=\mathbb N\cap [P+1,...,P+k]\hskip 0,1 cm $ by $\hskip 0,1 cm -P,\hskip 0,1 cm$ and leaves fixed all $\hskip 0,1 cm n\in \mathbb Z, \hskip 0,1 cm$ which are not in the disjoint union  $\hskip 0,1 cm  I\cup (P+I)\hskip 0,1 cm$, that is, which is   defined by
\begin{eqnarray}
\sigma(j)=P+j, \hskip 0,2 cm \textit{and}  \hskip 0,2 cm  
 \sigma(P+j)=j,\hskip 0,2 cm \textit{if} \hskip 0,2 cm 1\le j\le k\\
\sigma(n)=n\hskip 0,2 cm \textit{if} \hskip 0,2 cm n\notin\{1,...,k\}\cup \{P+j:j=1,...,k\},
\end{eqnarray}
so that
\begin{eqnarray}
 \nonumber T_\sigma(\omega)_n=\omega_n , \hskip 0,1 cm \textit{if}\hskip 0,2 cm n\notin\{1,...,k\}\cup \{P+j:j=1,...,k\}\\
T_\sigma(\omega)_j=\omega_{P+j}\hskip 0,5 cm  \textit{and} \hskip 0,5 cm  T_\sigma(\omega)_{P+j}=\omega_j  \hskip 0,5 cm \textit{for }\hskip 0,5 cm 1\le j\le k .
\end{eqnarray}

In this particular case, where $\sigma=\sigma_{P,k}$, which will be used  in the proof of Theorem 1,  we shall denote $\Gamma_\sigma$ by $\Gamma_{P,k}$, or, when $\hskip 0,1 cm k\hskip 0,1 cm$  is fixed, simply by $\hskip 0,1 cm \Gamma_P.\hskip 0,1 cm$ We have then the following equality
\begin{eqnarray*}
\Gamma_P=\bigcup_{(a,b),(c,d)\in {\cal Q}} C_k(a,b,c,d)\cap S^{-P}D_k(a,b,c,d).\hskip 1 cm (*)
\end{eqnarray*}
We note that, by $(19)$ and $(22)$,\begin{eqnarray} \forall  k>k_0, \forall  P>k\vee P_0(k), \hskip 0,4 cm \textit{the involution}\hskip 0,4 cm  \sigma_{P,k} \hskip 0,4 cm \textit{ belongs to} \hskip 0,4 cm H_+.\hskip 0,2 cm
\end{eqnarray} 
 To prove the equality  $(*)$,  let $x\in \Sigma$ and consider $T_\sigma x$, where $\sigma=\sigma_{P,k}$. Write them as an infinite strings
\begin{eqnarray*}
...x_{-1}x_0x_1x_2...x_kx_{k+1}...x_{P-1}x_Px_{P+1}...x_{P+k-1}x_{P+k}x_{P+k+1}...\\
...x_{-1}x_0x_{P+1}x_{P+2}...x_{P+k}x_{k+1}...x_{P-1}x_Px_1...x_{k-1}x_kx_{P+k+1}...
\end{eqnarray*}
the first line representing $x$, and the second corresponds to $T_\sigma x$. 
Then we have the following implications 
 \begin{eqnarray*}
x\in \Gamma_P\iff x\in \Sigma, T_\sigma x\in \Sigma\iff \\ x\in \Sigma, x_{P+1}\in R(x_0), x_{P+k}\in L(x_{k+1}), x_P\in L(x_1), x_{P+k+1}\in R(x_k)\\
\iff \exists a,b,c,d, x_0=a,x_1=b,x_k=c,x_{k+1}=d,\\ (S^Px)_1\in Ra, (S^Px)_k\in Ld, (S^Px)_0\in Lb, (S^Px)_{k+1}\in Rc\\
\iff \exists a,b,c,d, x\in C_k(a,b,c,d), S^px\in D_k(a,b,c,d)\\\iff  \exists a,b,c,d, x\in C_k(a,b,c,d), x\in S^{-P}\in D_k(a,b,c,d)
\end{eqnarray*}
which prove the equality $(*)$.\\
Let ${\cal T}$ be the two-sided  tail sigma field
 \begin{eqnarray*}
{\cal T}:=\bigcap_{n\ge 1}({\cal B}_{\le -n}\vee {\cal B}_{\ge n}).
\end{eqnarray*}

\section{The main result}
{\bf{Definition 2}}\\
\textit{ A probability measure $\mu$ on $\Sigma$ is said quasi invariant  if for any involution $\sigma$, the restrictions to $\Gamma_\sigma:=\Sigma\cap T_\sigma^{-1}\Sigma,$ of the two measures $\mu$ and $\mu\circ T_\sigma^{-1}$,  are equivalent.}\\ 

\textit{In this case, we still, as  in the full shift case,    let   $\phi_\sigma$ denote the Radon Nikodym derivative of $\hskip 0,1 cm \mu\circ T_\sigma^{-1}\hskip 0,1 cm $ restricted to $\hskip 0,1 cm \Gamma_\sigma\hskip 0,1 cm $ with respect to the restriction of $\hskip 0,1 cm \mu\hskip 0,1 cm $ to the same set.}\\
In the particular case where $\hskip 0,1 cm \sigma=\sigma_{P,k},\hskip 0,1 cm$ $\hskip 0,1 cm \phi_\sigma\hskip 0,1 cm $ will be denoted $\hskip 0,1 cm \phi_{\sigma_{P,k}}$. \\
{\bf{Theorem 1 }}\\
\textit{ Let $\hskip 0,1 cm (\Sigma,S)\hskip 0,1 cm$ be topologically mixing subshift. Let $\hskip 0,1 cm \mu\hskip 0,1 cm $ be a  shift invariant probability measure on $\hskip 0,1 cm \Sigma,\hskip 0,1 cm$   such that the system $\hskip 0,1 cm (\Sigma,S,\mu)\hskip 0,1 cm $ is ergodic. Suppose that  $\hskip 0,1 cm \mu\hskip 0,1 cm $ is quasi invariant. Then $\hskip 0,1 cm (\Sigma,S,\mu)\hskip 0,1 cm $ is isomorphic to a Bernoulli shift.} \\

{\bf{Proof}}:  The natural numbers $P$ and $k$ in this proof are as in $(23)$. Let $R\in \mathbb N$ be such that $R+k<P$.  For any  positive $n\in \mathbb N$, such that $n\le R$, let   $\hskip 0,1 cm V_1,...,V_n\hskip 0,1 cm $ be subsets of $\hskip 0,1 cm K.\hskip 0,1 cm$
Using the notation  as in (5) and (6), and setting
   \begin{eqnarray}
E=\Pi(A_{-L},...,A_0):=\{\omega\in \Sigma: \omega_{-L}\in A_{-L},...,\omega_0\in A_0\} ,
\end{eqnarray}

then from the equality $\hskip 0,1 cm T_\sigma^{-1}\Gamma_P=\Gamma_P,\hskip 0,1 cm$ from the definition of $ \hskip 0,1 cm \sigma=\sigma_{P,k}\hskip 0,1 cm $ given in $(20)$ and $(21)$, and from the definition of $\hskip 0,1 cm T_\sigma,\hskip 0,1 cm $ we obtain the  following equality
 for all $\hskip 0,1 cm T\in{\cal T},\hskip 0,1 cm$ and $\hskip 0,1 cm E\hskip 0,1 cm $ as in $(24)$,
\begin{eqnarray*}
\mu(T_\sigma^{-1}(\Gamma_P\cap E\cap F(A_1,...,A_k)\cap S^{-k} F(V_1,...,V_n)\cap S^{-P}F(B_1,...,B_k)\cap T)=\\\mu(\Gamma_P\cap E\cap F(B_1,...,B_k)\cap S^{-k}F(V_1,...,V_n)\cap S^{-P}F(A_1,...,A_k)\cap T)
\end{eqnarray*}
because of  the following two equalities
\begin{eqnarray*}
T_\sigma^{-1}( E\cap F(A_1,...,A_k)\cap S^{-k} F(C_1,...,C_M)\cap S^{-P}F(B_1,...,B_k)=\\ E\cap F(B_1,...,B_k)\cap S^{-k}F(C_1,...,C_M)\cap S^{-P}F(A_1,...,A_k)
\end{eqnarray*}
and  
\begin{eqnarray*}
T_\sigma^{-1}T=T.
\end{eqnarray*}
Recall that
\begin{eqnarray}
\Gamma_P=\bigcup_{(a,b),(c,d)\in {\cal Q}} C_k(a,b,c,d)\cap S^{-P}D_k(a,b,c,d),
\end{eqnarray}
 and let us  find, modulo $\mu$ zero set,  the following subset
\begin{eqnarray*}
R=R_{a,b,c,d}:=T_\sigma^{-1}(C_k(a,b,c,d)\cap S^{-P}D_k(a,b,c,d)).
\end{eqnarray*}
First 
\begin{eqnarray*} T_\sigma x\in C_k(a,b,c,d)\iff x_0=a, x_{P+1}=b,x_{P+k}=c, x_{k+1}=d\\
\iff x_0=a,x_{k+1}=d, x_{P+1}=b,x_{P+k}=c
\end{eqnarray*}
and, if we set $y=T_\sigma x$, 
\begin{eqnarray*}
T_\sigma x\in S^{-P}D_k(a,b,c,d)\iff S^PT_\sigma x\in D_k(a,b,c,d)\iff \\
(S^Py)_0\in Lb, (S^Py)_1\in Ra, (S^Py)_k\in Ld,(S^Py)_{k+1}\in Rc\\
\iff y_P\in Lb, y_{P+1}\in Ra, y_{P+k}\in Ld, y_{P+k+1}\in Rc\\
\iff x_P\in Lb, x_1\in Ra, x_k\in Ld, x_{P+k+1}\in Rc\\
\iff x_1\in Ra, x_k\in Ld, x_P\in Lb, x_{P+k+1}\in Rc
\end{eqnarray*}
So that
\begin{eqnarray*}
x\in R\iff x_0=a,x_1\in Ra, x_k\in Ld, x_{k+1}=d,x_P\in Lb, x_{P+1}=b,\\ x_{P+k}=c,x_{P+k+1}\in Rc
\end{eqnarray*}
so that modulo $\mu$ null set
\begin{eqnarray*}
x\in R\iff x_0=a,  x_{k+1}=d, x_{P+1}=b,\\ x_{P+k}=c
\end{eqnarray*}
that is
\begin{eqnarray*}
R=(x_0=a,x_{k+1}=d)\cap S^{-P}(x_1=b,x_k=c)
\end{eqnarray*}
It follows
\begin{eqnarray*}
R\cap T_\sigma^{-1}( E\cap F(A_1,...,A_k)\cap S^{-k} F(V_1,...,V_n)\cap S^{-P}F(B_1,...,B_k)=\\ R\cap E\cap F(B_1,...,B_k)\cap S^{-k}F(V_1,...,V_n)\cap S^{-P}F(A_1,...,A_k)
\end{eqnarray*}
so that
\begin{eqnarray*}
\int_{T\cap C_k(a,b,c,d)\cap S^{-P}D_k(a,b,c,d)\cap  E\cap F(A_1,...,A_k)\cap S^{-k} F(V_1,...,V_n)\cap S^{-P}F(B_1,...,B_k)}\phi_\sigma d\mu=\\
\mu(T\cap R\cap E\cap F(B_1,...,B_k)\cap S^{-k}F(V_1,...,V_n)\cap S^{-P}F(A_1,...,A_k))
\end{eqnarray*}
which we write as
\begin{eqnarray}
\nonumber \int_{T\cap E\cap F(A_1,...,A_k)\cap S^{-k}F(V_1,...,V_n)} 1_{C_k(a,b,c,d)} [1_{D_k(a,b,c,d)}1_{F(B_1,...,B_k)}]\circ S^P\phi_\sigma d\mu\\
=\int_{T\cap E\cap F(B_1,...,B_k)\cap S^{-k}F(V_1,...,V_n)}1_{(x_0=a,x_{k+1}=d)}[1_{(x_1=b,x_k=c)}1_{F(A_1,...,A_k)}]\circ S^P d\mu.
\end{eqnarray}
Define two sequences $\delta^Q_{(A_1,...,A_k)}(b,c)$ and $\xi^Q_{B_1,...,B_k}(a,b,c,d)$, $Q\ge R+k+1$, by 
\begin{eqnarray*}
\delta^Q_{(A_1,...,A_k)}(b,c):={1\over Q}\sum_{P=R+k+1}^Q [1_{(x_1=b,x_k=c)}1_{F(A_1,...,A_k)}]\circ S^P 
\end{eqnarray*}
and,  recalling that $\sigma=\sigma_{P,k}$,
\begin{eqnarray*}
\xi^Q_{B_1,...,B_k}(a,b,c,d):={1\over Q}\sum_{P=R+k+1}^Q [1_{D_k(a,b,c,d)}1_{F(B_1,...,B_k)}]\circ S^P\phi_\sigma
\end{eqnarray*}
so that, by $(26)$, for any $Q\ge R+k+1$, 
\begin{eqnarray}
\nonumber \int_{T\cap E\cap F(A_1,...,A_k)\cap S^{-k}F(V_1,...,V_n)} 1_{C_k(a,b,c,d)} \xi^Q_{B_1,...,B_k}(a,b,c,d)d\mu\\
=\int_{T\cap E\cap F(B_1,...,B_k)\cap S^{-k}F(V_1,...,V_n)}1_{(x_0=a,x_{k+1}=d)}\delta^Q_{(A_1,...,A_k)}(b,c)d\mu.
\end{eqnarray}
The sequence $ \hskip 0,1 cm (\xi^Q_{B_1,...,B_k}(a,b,c,d)_{Q\ge R+k+1}\hskip 0,1 cm$ 
is bounded in $L^\infty(\mu)^*$. By Alaoglu-Bourbaki theorem ([3] Theorem 2., p. 424), it then admits   at least one weak-star cluster point.  Because, by the mean ergodic Theorem, the sequence  $\hskip 0,1 cm\delta^Q_{(A_1,...,A_k)}(b,c)\hskip 0,1 cm $ converges in $L^2(\mu)$ norm, the eququality $(27)$ implies that for any weak star cluster point $\hskip 0,1 cm \eta_{B_1,...,B_k}(a,b,c,d)\hskip 0,1 cm $  of the sequence $\hskip 0,1 cm (\xi^Q_{B_1,...,B_k}(a,b,c,d))_{Q\ge R+k+1}\hskip 0,1 cm $, we have
\begin{eqnarray*}
\eta_{B_1,...,B_k}(a,b,c,d)(C_k(a,b,c,d)\cap T\cap E\cap F(A_1,...,A_k)\cap S^{-k}F(V_1,...,V_n))\\
=\mu((x_1=b,x_k=c)\cap F(A_1,...,A_k))\\\times \mu((x_0=a,x_{k+1}=d)\cap T\cap E\cap F(B_1,...,B_k)\cap S^{-k}F(V_1,...,V_n))
\end{eqnarray*}
In the particular case where $B_1=...=B_k=K$, we denote $\eta_{B_1,...,B_k}(a,b,c,d)$ by $\eta(a,b,c,d)$. Also the dependence on $k$ is specified by $$\eta^k=\eta^k(a,b,c,d)=\eta(a,b,c,d).$$
Then, for $E=E_{-}\cap (x_0=a)$, where $E_-$ is of the form 
\begin{eqnarray*}
E_-=\Pi(A_{-L},...,A_{-1},K):=\{\omega\in \Sigma: \omega_{-L}\in A_{-L},...,\omega_{-1}\in A_{-1}\}
\end{eqnarray*}
we have
\begin{eqnarray*}
\eta^k[(x_0=a,x_1=b, x_k=c,x_{k+1}=d)\cap T\cap  E_{-}\cap (x_1\in A_1,...,x_k\in A_k,x_{k+1}\in V_1,..., x_{k+n}\in V_n)]=\\
\mu[(x_1=b,x_k=c, x_1\in A_1,...,x_k\in A_k)]\times\\ \mu[(x_0=a,x_{k+1}=d)\cap T\cap E_{-}\cap (x_0=a)\cap (x_{k+1}\in V_1,...,x_{k+n}\in V_n)]  
\end{eqnarray*}
which, when $b\in A_1, c\in A_k$ and $d\in V_1$,  we rewrite as
\begin{eqnarray*}
\eta^k[ T\cap E_{-}\cap (x_0=a,x_1=b, x_2\in A_2,...,x_{k-1}\in A_{k-1},  x_k=c,x_{k+1}=d,x_{k+2}\in V_2,..., x_{k+n}\in V_n)]=\\
\mu(x_1=b, x_2\in A_2,...,x_{k-1}\in A_{k-1}, x_k=c)\times \mu(T\cap E_{-}\cap(x_0=a,x_{k+1}=d, x_{k+2}\in V_2,...,x_{k+n}\in V_n)
\end{eqnarray*}
In particular, for any $l$ with $2\le l\le k-1$, 
\begin{eqnarray*}
\eta^k[ T\cap E_{-}\cap (x_0=a,x_1=b, x_2\in A_2,...,x_l\in A_l,  x_k=c,x_{k+1}=d,x_{k+2}\in V_2,..., x_{k+n}\in V_n)]=\\
\mu(x_1=b, x_2\in A_2,...,x_l\in A_l, x_k=c)\times \mu(T\cap E_{-}\cap(x_0=a,x_{k+1}=d, x_{k+2}\in V_2,...,x_{k+n}\in V_n)
\end{eqnarray*}
that is
\begin{eqnarray*}
\eta^k[ T\cap E_{-}\cap (x_0=a,x_1=b, x_2\in A_2,...,x_l\in A_l)\cap S^{-k}( x_0=c,x_1=d,x_2\in V_2,..., x_n\in V_n)]=\\
\mu((x_1=b, x_2\in A_2,...,x_l\in A_l)\cap S^{-k}( x_0=c))\\\times \mu(T\cap E_{-}\cap(x_0=a)\cap S^{-k}(x_1=d, x_2\in V_2,...,x_n\in V_n))
\end{eqnarray*}
Let 
\begin{eqnarray*}
\alpha_{a,b}^k=\alpha^k:=\sum_{c,d} \eta^k(a,b,c,d).
\end{eqnarray*}
Then, in particular for $V_2=...=V_n=K$, we obtain
\begin{eqnarray*}
\alpha_{a,b}^k(T\cap E_{-}\cap (x_0=a,x_1=b,x_2\in A_2,...,x_l\in A_l)=\\\mu( x_1=b,x_2\in A_2,...,x_l\in A_l)\times \mu(T\cap E_{-}\cap (x_0=a)).
\end{eqnarray*}
Clearly the sequence $\hskip 0,1 cm (\eta^k(a,b,c,d))_{k\ge 3}\hskip 0,1 cm $  is bounded in $L^{\infty *}(\mu). $ The same holds then for the sequence $\hskip 0,1 cm (\alpha^k_{a,b})_{k\ge 3}.\hskip 0,1 cm$  
It follows that every weak star cluster point $\hskip 0,1 cm \alpha_{a,b}\in L^{\infty*}(\mu)\hskip 0,1 cm$ of the sequence $\hskip 0,1 cm (\alpha_{a,b}^k)_{k\ge 3},\hskip 0,1 cm$ satisfies: for any $l\ge 2$, 
\begin{eqnarray*}
\alpha_{a,b}(T\cap E_{-}\cap (x_0=a,x_1=b,x_2\in A_2,...,x_l\in A_l)=\\ \mu(T\cap E_{-}\cap (x_0=a))\times \mu( x_1=b,x_2\in A_2,...,x_l\in A_l)
\end{eqnarray*}
so that $\hskip 0,1 cm \alpha_{a,b}\hskip 0,1 cm $ extends uniquely to a countably additive $\hskip 0,1 cm \tilde{\alpha}_{a,b} \hskip 0,1 cm $ to the algebra of the sets of the form $\hskip 0,1 cm  E\cap (x_0=a,x_1=b)\cap F,\hskip 0,1 cm$ where $\hskip 0,1 cm E\in \bigvee_{j\le -1} {\cal B}_j,\hskip 0,1 cm $  $\hskip 0,1 cm F\in \bigvee_{j\ge 2} {\cal B}_j,\hskip 0,1 cm $ and
\begin{eqnarray}
\tilde{\alpha}_{a,b}( T\cap E\cap (x_0=a,x_1=b)\cap F)=\mu(T\cap E\cap (x_0=a))\times \mu((x_1=b)\cap F).
\end{eqnarray}
 
 Also this last equality still holds for $\hskip 0,1 cm E\hskip 0,1 cm $ in the $\mu$ completion of $\hskip 0,1 cm \bigvee_{j\le -1} {\cal B}_j\hskip 0,1 cm $ and $\hskip 0,1 cm F\hskip 0,1 cm $ in  the $\mu$ completion of $\hskip 0,1 cm \bigvee_{j\ge 2} {\cal B}_j.\hskip 0,1 cm$ 
[In particular it still holds for all $\hskip 0,1 cm E,F\in \Pi_S,\hskip 0,1 cm$ where $\hskip 0,1 cm \Pi_S\hskip 0,1 cm$ is the Pinsker sigma algebra of the system $\hskip 0,1 cm(\Omega,S,\mu),\hskip 0,1 cm $so that we have\\
{\bf{Proposition 1 }}\\
\textit{The Pinsker sigma algebra of the system $(\Omega,S,\mu)$ is trivial.}\\
{\bf{Proof}}\\ In fact if $R\in \Pi_S$ with $\mu(R)>0$,  then
\begin{eqnarray*}
\tilde{\alpha}( R\cap (x_0=a,x_1=b)\cap R)=\mu(R\cap (x_0=a))\times \mu((x_1=b)\cap R)\\\textit{and also} \hskip 0,5 cm 
\tilde{\alpha}( R\cap (x_0=a,x_1=b)\cap R)=\mu(R\cap (x_0=a))\times \mu((x_1=b))\\
\textit{imply}\hskip 0,5 cm \mu(R\cap (x_0=a))\times \mu((x_1=b)\cap R)=\mu(R\cap (x_0=a))\times \mu((x_1=b)).
\end{eqnarray*}
Then, for $a$ with $ \mu(R\cap (x_0=a))\ne 0$, we have 
\begin{eqnarray*}
\mu((x_1=b)\cap R)=\mu((x_1=b)), \forall b\\
\textit{thus}\hskip 0,5 cm \mu(R)=1.
\end{eqnarray*}]

Now, we continue the proof of Theorem 1. Let 
\begin{eqnarray*}
\nu:=\sum_{a,b} \tilde{\alpha}_{a,b}.
\end{eqnarray*}
Then from $(28)$ we get 
\begin{eqnarray}
\nu(T\cap E\cap F)=\mu(T\cap E)\times \mu(F), \hskip 0,5 cm \forall T\in{\cal T}, \hskip 0,1 cm  E\in {\cal B}_{\le 0}, \hskip 0,1 cm  F\in {\cal B}_{\ge 1}.
\end{eqnarray}
Let ${\cal G}:={\cal B}_{\le 0}\times{\cal B}_{\ge 1}$ and ${\cal L}:=({\cal T}\vee {\cal B}_{\le 0})\times {\cal B}_{\ge 1}$. Then $(29)$ shows that $\nu$ is countably additive on ${\cal L}$. Thus $\nu$ extends uniquely to a countably additive measure $\nu_1$ to the sigma algebra $\hskip 0,1 cm \sigma({\cal L})\hskip 0,1 cm$ generated by $\hskip 0,1 cm {\cal L}.\hskip 0,1 cm$ Since $\hskip 0,1 cm {\cal G}\subset {\cal L},\hskip 0,1 cm$  $\hskip 0,1 cm \nu\hskip 0,1 cm $ extends uniquely to a countably additive measure  $\hskip 0,1 cm \nu_2\hskip 0,1 cm $ on $\hskip 0,1 cm \sigma({\cal G})\hskip 0,1 cm$ also. 
Because $\hskip 0,1 cm \sigma({\cal G})=\sigma({\cal L})={\cal B},\hskip 0,1 cm$ and $\hskip 0,1 cm \nu_1\hskip 0,1 cm $ is a countably additive extension of $\hskip 0,1 cm \nu\hskip 0,1 cm$ from $\hskip 0,1 cm {\cal G}\hskip 0,1 cm$ to $\hskip 0,1 cm \sigma({\cal G}),\hskip 0,1 cm$ we have $\hskip 0,1 cm \nu_1=\nu_2=:\tilde{\nu}.\hskip 0,1 cm$ Then in particular $\hskip 0,1 cm \tilde{\nu}\hskip 0,1 cm$ is countably additive on $\hskip 0,1 cm {\cal B}\hskip 0,1 cm$ and satisfies the two following equalities
\begin{eqnarray}
\tilde{\nu}( E\cap F)=\mu( E)\times \mu(F), \hskip 0,5 cm  \hskip 0,1 cm  \forall E\in {\cal B}_{\le 0}, \hskip 0,1 cm  F\in {\cal B}_{\ge 1}.\\
\textit{and} \hskip 1 cm \tilde{\nu}(T)=\mu(T),\hskip 0,5 cm \forall T\in{\cal T} 
\end{eqnarray}
which mean that $\mu$ is "faiblement de Bernoulli". Since "faiblement de Bernoulli" is equivalent to  weak Bernoulli ( [7], Proposition 2) and also, a system which is weak Bernoulli is isomorphic to a Bernoulli system [4],  the proof is complete.

\newpage

\section{Application}

\subsection {Quasi-invariance of Gibbs measures}If $\Lambda$ is an $n\times n$ matrix of zeros and ones, let $K:=\{0,1,...,n-1\}$ and 
\begin{eqnarray*}
\Sigma_\Lambda:=\{x\in K^\mathbb Z: \Lambda_{x_i,x_{i+1}}=1, \forall i\in \mathbb Z\}.
\end{eqnarray*}
We asume that $\forall j\in K$, there exists $x\in \Sigma_\Lambda$ such that $x_0=j$. \\

{\bf{Definition 3}}\\  Let  $\phi:\Sigma_\Lambda\rightarrow \mathbb R$ be  continuous.  \textit{A Gibbs measure for $\phi$ is a shift invariant probability measure $\mu_\phi$ on $\Sigma_\Lambda$ for which one can find constants $c_1>0, c_2>0$ and p such that
\begin{eqnarray}
c_1\le {\mu_\phi(\{y\in \Sigma_\Lambda: y_j=x_j, j=0,...,m\})\over 
\textit{exp}(-pm+\sum_{j=0}^{m-1} \phi(S^jx))}\le c_2,
\end{eqnarray}
 for every $\hskip 0,1 cm x\in \Sigma_\Lambda\hskip 0,1 cm $ and $\hskip 0,1 cm m\ge 0$.}\\
 
 Recall that\textit{ the system $\hskip 0,1 cm (\Sigma_\Lambda,S)\hskip 0,1 cm $ is topologically mixing if for some $\hskip 0,1 cm M,\hskip 0,1 cm$ $\hskip 0,1 cm \Lambda^M_{i,j}>0,\hskip 0,1 cm$ for all $\hskip 0,1 cm i,j$}.\\
 Let, as  in [1],   $\hskip 0,1 cm {\cal{F}}_\Lambda\hskip 0,1 cm $ be the set of $\hskip 0,1 cm \phi\hskip 0,1 cm$ which satisfies  
 \begin{eqnarray*}
\textit{var}_k\phi \le b\alpha^k, \forall k\ge 0,
\end{eqnarray*}
 for some $\hskip 0,1 cm \alpha\in ]0,1[\hskip 0,1 cm $ and $\hskip 0,1 cm b,$ 
  where   $\hskip 0,1 cm \textit{var}_k\phi\hskip 0,1 cm $ is defined by
\begin{eqnarray*}
\textit{var}_k\phi:=\sup\{\mid \phi(x)-\phi(y)\mid :x,y\in \Sigma_\Lambda, x_j=y_j,\forall \mid j\mid\le k\}.
\end{eqnarray*}
Clearly, $\hskip 0,1 cm{\cal{F}}_\Lambda\hskip 0,1 cm $ contains every $\hskip 0,1 cm \phi\hskip 0,1 cm $ which depends only on a finite number of coordinates.\\ 

 {\bf{Definition 4}}:\\ \textit{Two functions $\phi,\psi\in C(\Sigma_\Lambda)$ are homologous with respect to the shift $S$, if there is a $u\in C(\Sigma_\Lambda)$ such that \begin{eqnarray*}
\psi(x)=\phi(x)+u(x)-u\circ S(x), \hskip 0,3 cm \forall x\in \Sigma_\Lambda.
\end{eqnarray*}}
We recall the following results from [1]:\\
{\bf{Theorem A}}:\\
\textit{Suppose  $(\Sigma_\Lambda,S)$  topologically mixing and let $\phi\in {\cal{F}}_\Lambda.$  Then\\
 (i)   there exists a unique Gibbs measure $\mu_\phi$  for $\phi$.\\
 (ii) If $\psi$ is cohomologous to $\phi$ then $\mu_\psi=\mu_\phi$.\\
 (iii)  $\phi$ is  cohomologous to some  $\psi\in {\cal{F}}_\Lambda$ with $\psi(x)=\psi(y)$ whenever $x_j=y_j$ for all $j\ge 0$. }\\
 (iv)   $(\Sigma_\Lambda,S,\mu_\phi)$ is isomorphic to a Bernoulli system.\\
  
We prove the following proposition, from which the statement (iv) in Theorem A, follows as a corollary of Theorem 1:\\
{\bf{Proposition 2}}\\{\textit{For every potential $\hskip 0,1 cm \phi\in {\cal {F}}_\Lambda\hskip 0,1 cm $ the Gibbs measure $\hskip 0,1 cm \mu_\phi\hskip 0,1 cm$ is quasi-invariant by any involution which moves only a finite number of coordinaotes. }\\

{\bf{Proof}} 
 Denote  $\hskip 0,1 cm \mu_\phi\hskip 0,1 cm $ by $\hskip 0,1 cm \mu$. Observe first that the  inequalities   (32) in the definition can be written as 
\begin{eqnarray}
\mu(\{y\in \Sigma_\Lambda: y_j=x_j, j=0,...,m\})=a_m(x) \textit{exp}(-pm+\sum_{j=0}^{m-1} \phi(S^jx)),
\end{eqnarray}
with
\begin{eqnarray*}
c_1\le a_m(x)\le c_2, \end{eqnarray*}
and, in particular,  implie that the measure of any cylinder is non nul.\\ 
  Let 
  \begin{eqnarray*}H_\Lambda :=\{\sigma\in H: \sigma \hskip 0,2 cm \textit{is an involution},\hskip 0,2 cm T_\sigma \Sigma_\Lambda \cap \Sigma_\Lambda \ne\emptyset\}
  \end{eqnarray*}.  Then 
  \begin{eqnarray*}
\sigma\in H_\Lambda\iff \sigma^{-1}\in H_\Lambda,
\end{eqnarray*}
because of the equality
\begin{eqnarray*}
T_{\sigma^{-1}}(\Sigma_\Lambda \cap T_\sigma \Sigma_\Lambda )=T_{\sigma^{-1}}\Sigma_\Lambda  \cap \Sigma_\Lambda .
\end{eqnarray*}

Let  $\sigma\in H_\Lambda,\hskip 0,1 cm$  $\hskip 0,1 cm \tau=\sigma^{-1}\hskip 0,1 cm $ and  $\hskip 0,1 cm M(\tau)\hskip 0,1 cm $ be the smallest natural number such that $$\mid j\mid \ge M(\tau)\Rightarrow \tau(j)=j.$$ We note then that, for any $\hskip 0,1 cm N\ge M(\tau), \hskip 0,1 cm$  the restriction of $\hskip 0,1 cm \tau\hskip 0,1 cm$ to the set $\hskip 0,1 cm\{-N,...,N\}\hskip 0,1 cm$ is a permutation of this set.    Let $\hskip 0,1 cm x\in \Sigma_\Lambda\cap T_\sigma^{-1}\Sigma_\Lambda,\hskip 0,1 cm$ and for all $\hskip 0,1 cm N\ge M(\tau),\hskip 0,1 cm$ \begin{eqnarray*}
C_N(x)=C:=\{y\in \Sigma_\Lambda: y_j=x_j, j=-N,...,N\}.
\end{eqnarray*}
Then, by shift-invariance of $\hskip 0,1 cm \mu$,  
\begin{eqnarray*}
\mu(C)=
=\mu(\{y\in \Sigma_\Lambda:y_0=x_{-N},y_1=x_{-N+1},...,y_{2N}=x_N\},
\end{eqnarray*}
hence by $(33)$, we obtain
\begin{eqnarray}
\nonumber \mu(C_N(x))=\mu(C)=a_{2N}(S^{-N}x)\textit{exp}(-p\times 2N+\sum_{j=0}^{2N-1}\phi(S^jS^{-N}x)
\end{eqnarray}
that is 
\begin{eqnarray}
\mu(C_N(x))= a_{2N}(S^{-N}x)\textit{exp} (-2Np+\sum_{q=-N}^{N-1} \phi(S^q x))
\end{eqnarray}

and also, since, as noted before, $\hskip 0,1 cm \sigma\hskip 0,1 cm$ restricted to $\hskip 0,1 cm \{-N,...,N\}\hskip 0,1 cm $ is a permutation of $\hskip 0,1 cm \{-N,...,N\},\hskip 0,1 cm $ we have  
\begin{eqnarray*}
T_\sigma^{-1}C=\{y\in \Sigma_\Lambda:y_{\sigma(j)}=x_j,j=-N,...,j=N\}\\
=\{y\in \Sigma_\Lambda:y_{j}=x_{\tau(j)},j=-N,...,j=N\},
\end{eqnarray*}
where $\hskip 0,1 cm \tau :=\sigma^{-1},\hskip 0,1 cm $ and because  $\hskip 0,1 cm x_{\tau (j)}=(T_\tau x)(j),\hskip 0,1 cm$  for all $\hskip 0,1 cm j,\hskip 0,1 cm$  we obtain
\begin{eqnarray*}
T_\sigma^{-1}C=C_N(T_\tau x),
\end{eqnarray*}
 so by $(34)$, 
\begin{eqnarray*}
\mu(T_\sigma^{-1}C)=a_{2N}( S^{-N}T_\tau x)\textit{exp}(-2Np+\sum_{q=-N}^{N-1} \phi(S^q T_\tau x)).
\end{eqnarray*}

It follows that
\begin{eqnarray}
{\mu(T_\sigma^{-1}C)\over \mu(C)}={a_{2N}(  S^{-N} T_\tau x)\over a_{2N(S^{-N}x)}}\times \textit{exp}(\sum_{q=-N}^{N-1} ( \phi(S^q  T_\tau x)-\phi(S^q x))
).\end{eqnarray}
Set, $\hskip 0,1 cm \tau\hskip 0,1 cm$ being fixed,  
\begin{eqnarray*}
G_{N,\tau}(x)=G_N(x):=\sum_{q=-N}^{N-1} ( \phi(S^q  T_\tau x)-\phi(S^q x),
\end{eqnarray*}
so that $(35)$ becomes 
\begin{eqnarray}
{\mu(T_\sigma^{-1}C)\over \mu(C)}={a_{2N}(  S^{-N} T_\tau x)\over a_{2N(S^{-N}x)}}\times \textit{exp}( G_N(x)).
\end{eqnarray}

 But, according to Theorem A, we can suppose that $\hskip 0,1 cm \phi\hskip 0,1 cm $ depends only on the non negative coordinates $\hskip 0,1 cm x_0,x_1,...$
\begin{eqnarray}
\phi(x)=\phi_1(x)=\phi_1(x_0,x_1,...).
\end{eqnarray}
Let  $\hskip 0,1 cm b>0\hskip 0,1 cm $ and $\hskip 0,1 cm \alpha\in]0,1[\hskip 0,1 cm $ such that  $\hskip 0,1 cm  \textit{var}_k \phi_1 \le b\alpha^k,\hskip 0,1 cm  \forall k\ge 0.\hskip 0,1 cm$
 Then
\begin{eqnarray*}
G_N(x)=\sum_{q=-N}^{N-1} (\phi_1 (S^qT_\tau x)-\phi_1(S^q x))=\sum_{q=-N}^{N-1} (\phi_1(x_{\tau(q)},x_{\tau(q+1)},...)- \phi_1(x_q,x_{q+1},...,)),
\end{eqnarray*}
so that if $\hskip 0,1 cm N\ge M(\tau)+1,\hskip 0,1 cm$ we can write
\begin{eqnarray*}
G_N(x)=H_N(x)+F_N(x),
\end{eqnarray*}
where 
\begin{eqnarray*}
H_N(x)=\sum_{q=-N}^{-M(\tau)} (\phi_1(x_{\tau(q)},x_{\tau(q+1)},...)- \phi_1(x_q,x_{q+1},...)),\\
F_N(x)=\sum_{q=-M(\tau)+1}^{N-1} (\phi_1(x_{\tau(q)},x_{\tau(q+1)},...)- \phi_1(x_q,x_{q+1},...)).
\end{eqnarray*}
Clearly, for any $N\ge M(\tau)$, we have 
\begin{eqnarray*}
F_N(x)=F_{M(\tau)}(x).
\end{eqnarray*}
Also, since for any $q\le -M(\tau)$, 
\begin{eqnarray*}
(S^qT_\tau x)_j=(S^qx)_j,\hskip 0,1 cm \forall  j, 0\le j\le -M(\tau)-q,
\end{eqnarray*}
we get
\begin{eqnarray*}
\mid \phi_1(S^qT_\tau x)-\phi_1(S^q x)\mid \le \textit{var}_{-M(\tau)-q} \phi_1\le b\alpha^{-M(\tau)-q},
\end{eqnarray*}
and then 
\begin{eqnarray*}
 \sum_{q=-N}^{-M(\tau)} \mid (\phi_1(x_{\tau(q)},x_{\tau(q+1)},...)- \phi_1(x_q,x_{q+1},...))\mid \\
 \le  \sum_{q=-N}^{-M(\tau)} \textit{var}_{-M(\tau)-q} \phi_1=\sum_{k=0}^{N-M(\tau)} \textit{var}_k  \phi_1 \le \sum_{k=0}^\infty  \textit{var}_k  \phi_1\le{b\over 1-\alpha}. \end{eqnarray*}
It follows that, the sequence $H_N(x)$ converges  and thus  $G_N(x)$ converges also to \begin{eqnarray}\lim_N G_N(x)=\sum_{q\in \mathbb Z} (\phi(S^qT_\tau x)-\phi(S^qx)),
\end{eqnarray} and
\begin{eqnarray*}
\mid G_N(x)\mid \le \mid F_{M(\tau)}(x)\mid +{b\over 1-\alpha},
\end{eqnarray*}
from which and because $F_{M(\tau)}$ is continuous on the compact space $\Sigma_\Lambda$, we conclude that for some constant $C$, 
\begin{eqnarray*}
\sup_{N\ge M(\tau)} \sup_{x} \mid G_N(x)\mid \le C,
\end{eqnarray*}

which, in view of $(36)$, gives
\begin{eqnarray}
{c_1\over c_2}\times \textit{exp}(-C)\le {\mu(T_\sigma^{-1}C)\over \mu(C)}\le {c_2\over c_1}\times \textit{exp}(C),
\end{eqnarray}
because of the two following inequalities
  \begin{eqnarray*}
{c_1\over c_2}\le  {a_{2N}(  S^{-N} T_\tau x)\over a_{2N(S^{-N}x)}}\le {c_2\over c_1}.
\end{eqnarray*}
Let $\hskip 0,1 cm  \alpha:={c_1\over c_2}\textit{exp}(-C)\hskip 0,1 cm $ and $\hskip 0,1 cm \beta:={c_2\over c_1}\textit{exp}(C),\hskip 0,1 cm $ so that $(39)$ reads as
\begin{eqnarray*}
\alpha \mu(C)\le \mu(T_\sigma^{-1}C)\le \beta \mu(C),
\end{eqnarray*}
for all cylinder $\hskip 0,1 cm C\hskip 0,1 cm $ in the algebra  generated by the coordinates in $\hskip 0,1 cm \{-N,...,N\},\hskip 0,1 cm $ and thus, by finite additivity, these inequalities still hold for any set $\hskip 0,1 cm E\hskip 0,1 cm$ in that algebra. Since   $\hskip 0,1 cm N\ge M(\tau)\hskip 0,1 cm$ is arbitrary,   
  the equivalence of $\hskip 0,1 cm \mu\hskip 0,1 cm $ and $\hskip 0,1 cm \mu\circ T_\sigma^{-1}\hskip 0,1 cm $ follows, and this proves the quasi-invariance of $\mu=\mu_\phi$.$\square$\\
  
{\bf{Remark }}\\
In [10] the authors define  a Gibbs measure $\mu$ of a map $\phi:\Sigma_\Lambda\rightarrow \mathbb R$, with summable variation, to be a probability measure $\hskip 0,1 cm \mu\hskip 0,1 cm $ on $\hskip 0,1 cm \Sigma_\Lambda\hskip 0,1 cm$  such,   in our setting,    that for any $\sigma\in H_\Lambda$,  the Radon-Nikodym derivative ${d\mu\circ T_\sigma\over d\mu}$ satisfies
\begin{eqnarray*}
\log( {d\mu\circ T_\sigma\over d\mu}(x))=\sum_{k\in \mathbb Z} (\phi(S^kT_\sigma x)-\phi(S^k x)),
\end{eqnarray*}
for $\mu$ almost all $x$. Then, according to this definition, a Gibbs measure is quasi-invariant.   They then give a proof  of the fact that, for every $\phi$, with summable variation, there exists a unique Gibbs measure $\mu_\phi$, such that the dynamical system $(\Sigma_\Lambda,S,\mu_\phi)$ is a $K$-system.   \\
 By (36) and  (38), recalling that $\tau=\sigma^{-1}$, we have 
 \begin{eqnarray*}
\log({d\mu\circ T_\sigma^{-1}\over d\mu}(x))=\lim_N \log({a_{2N}(  S^{-N} T_\tau x)\over a_{2N(S^{-N}x)}})+\sum_{q\in \mathbb Z} (\phi(S^qT_\tau x)-\phi(S^qx))\\
\textit{or equivalently}\hskip 1 cm \log({d\mu\circ T_\tau\over d\mu}(x))=\log(g_\tau(x))+\sum_{q\in \mathbb Z} (\phi(S^qT_\tau x)-\phi(S^qx))
\end{eqnarray*}
where \begin{eqnarray*}
{c_1\over c_2}\le g_\tau(x)\le {c_2\over c_1}
\end{eqnarray*}
so that  the Gibbs measures as in [10] are Gibbs measures according to the given Definition 3
 which is quoted from [1].\\

{\bf{Acknowledgments}}.  
Jean Paul Thouvenot posed to me the question  whether or not quasi-invariance implies complete positive entropy. He suggests many improvements and complements. I am greatly indebted to him. I am deeply grateful to him  for many useful discussions as well as for valuable comments which improve  the presentation of the paper,   and also for  bringing to my attention some  appropriate references.
\newpage \centerline{References}

1.   R. Bowen (1975): Equilibrium States and the Ergodic Theory of Anosov Diffeomorphisms. Lect. Notes in Math. 470, Springer.\\
2.  A. I. Bufetov (2018): Quasi-Symmetries of Determinantal  Point Processes. Ann. Prob., Vol. 46. N$\small {0}$. 2, 956-1003.\\
3.   N. Dunford and J. T. Schwartz : Linear operators, Part I. General Theory, Interscience publishers, New York.\\
4.   N. A. Friedman and  D. S. Ornstein (1970):  On isomorphism of weak Bernoulli transformations. Adv. in Math. 5, 365-394.\\
5.  D. Hamdan: Ergodic quasi-exchangeable stationary processes are isomorphic to Bernoulli processes. arxiv:1903;10804v1, 26 March 2019, arxiv:1903.10804v2, 1 July 2020. \\
6.  E. Hewitt and  L. J. Savage (1955): Symmetric measures on Cartesian products. Trans. Amer. Math. Soc. 80, 470-501.\\
7.  F. Ledrappier (1976):  Sur la Condition de Bernoulli Faible et ses Applications. Lect. Notes in Math., vol 532, 152-159, Springer-Verlag.\\
8.  D. Ornstein  (1970): Bernoulli shifts with the same entropy are isomorphic. Adv. Math. 4,337-352.\\
9.  D. Ornstein  (1971):  Two Bernoulli shifts with infinite entropy are isomorphic. Adv. Math. 5, 339-348.\\
10.  K. Petersen, K. schmidt (1997): Symmetric Gibbs Measures. Trans. Amer. Math. Soc. Vol. 349, Number 7, p. 2775-2811.\\

.\\

  {Sorbonne Universit\' e,  UMR 8001, Laboratoire de Probabilit\' es, Statistique et Mod\' elisation,  Bo\^ite courrier 158, 4 Place Jussieu, F-75752 Paris Cedex 05, France.\\
E-mail: doureid.hamdan@upmc.fr}
\newpage

\end{document}